\documentclass[leqno,11pt]{amsart}

\newtheorem{thm}{Theorem}[section]
\newtheorem{cor}[thm]{Corollary}
\newtheorem{lemma}[thm]{Lemma}

\theoremstyle{definition}

\theoremstyle{remark}
\newtheorem{remark}[thm]{Remark} 

\newcommand{\abs}[1]{\lvert#1\rvert}

\newcommand{\R}{\mathbb R}

\newcommand{\C}{\mathbb C}

\newcommand{\N}{\mathbb N}

\newcommand{\ie}{i.~e.,\ }

\begin{document}
%

	
\title[Positive QF] {Positive trigonometric Quadrature Formulas and quadrature
on the unit circle} \author[Franz Peherstorfer]{Franz Peherstorfer*$^1$}
\thanks{ *Supported by the Austrian Science Fund FWF, project no.\ P20413-N18}
\thanks{$^1$The last modifications and corrections of this manuscript were done by the author 
in July 2009 (submitted to arXiv by P. Yuditskii and I. Moale).} 
\begin{abstract}
	We give several descriptions of positive quadrature formulas which are exact for
	trigonometric - , respectively, Laurent polynomials of degree less or equal
	$n-1-m$, $0\leq m\leq n-1$. A complete and simple description is obtained with
	the help of orthogonal polynomials on the unit circle. In particular it is shown
	that the nodes polynomial can be generated by a simple recurrence
	relation. As a byproduct interlacing properties of zeros of para-orthogonal
	polynomials are obtained. Finally, asymptotics for the quadrature weights are
	presented.
\end{abstract}
\maketitle

%

\section{Introduction}
Let $l \in {\mathbb N}_0 := {\mathbb N} \cup \{ 0 \},$ $\gamma \in \{ 0,1/2 \}$,
and let $\mathcal{T}_{l,\gamma} = \{ \sum\limits_{k=0}^{l} a_k
\cos (k + \gamma)\varphi + b_k \sin (k + \gamma)\varphi : a_k, b_k \in \mathbb R \}.$
By  $\mathcal{T}_l :=  \mathcal{T}_{l,0}$ we denote the set of trigonometric
polynomials of degree less or equal $l.$ Further let $\sigma$ be a positive measure on
$[0,2\pi]$ normalized by $\frac{1}{2\pi}\int _0^{2\pi}d\sigma = 1$. We call a quadrature
formula (qf) of the form
\begin{equation}\label{1}
   \frac{1}{2 \pi} \int_{0}^{2 \pi} t(\varphi) d\sigma(\varphi) =
   \sum\limits_{s=1}^{n} \mu_s t(\varphi_s) + R_n(t)
\end{equation}
with $0 \leq \varphi_1 < \varphi_2 < ... < \varphi_n < 2 \pi$, $\mu_1,...,\mu_n
\in \mathbb R,$ a trigonometric $(n-1-m, n, d\sigma)$ qf if $R_n(t) = 0$
for $t \in {\mathcal T}_{n-1-m}$; $n-1-m$ is the so-called degree of exactness.
Note that for $n \in \mathbb N,$ $n = 2 (\tilde{n} + \gamma),$ $\gamma \in \{ 0, 1/2 \}$,
the trigonometric nodes polynomial $\prod_{s=1}^n \sin(\frac{\varphi - \varphi _s}{2}) = e^{-i (n \varphi + \sum\limits_{s=1}^{n} \varphi_s)/2}
\prod_{s=1}^{n}(e^{i \varphi} - e^{i \varphi_s})/2i$ is from $\mathcal{T}_{\tilde{n},\gamma}.$


Obviously qf \eqref{1} is a $(n-1-m, n, d\sigma)$ qf if and only if
\begin{equation}\label{1.5}
\frac{1}{2\pi}\int\limits_{0}^{2\pi} e^{-ik\varphi} d\sigma(\varphi) = \sum_{s=1}^n\mu_s
e^{-ik\varphi_s} \hspace{1cm} k=0,\ldots,\pm(n-1-m)
\end{equation}
which is a qf on the unit circle exact for Laurent polynomials $ \sum
_{-(n-1-m)}^{n-1-m}d_k z^k,$ $d_k \in \mathbb{C}$ of degree $\leq n-m-1$, which we
call $(n-1-m, n, d\sigma)$ Laurent qf also. In other words we are looking for a
characterization of quadrature formulas which integrate exactly the first
$n-1-m$ moments $c_1,\ldots,c_{n-1-m}$ of the Caratheodory function

\begin{equation}
F(z) = \frac{1}{2\pi}\int\limits_{0}^{2\pi}\frac{e^{i\varphi}+z}{e^{i\varphi}-z}d\sigma(\varphi)
= 1 + 2\sum_{k=1}^{\infty}c_kz^k.
\label{eq-f1}
\end{equation}

Because of convergence and stability reasons it is of upmost interest to have positive
$(n-1-m,n,d\sigma )$ qf, i.e., that all quadrature weights $\mu _s$ are positive. Such qf
will be studied in this paper. So far it is known that qf with highest possible
degree of exactness, that is, $(n-1,n,d\sigma )$ qf have positive weights.
These $(n-1,n,d\sigma )$ qf can be considered as the counterparts of the Gauss qf for algebraic
polynomials and are called Szeg\H{o} qf nowadays. The Szeg\H{o} qf can be described as follows
(see \cite[pp. 14-16]{Ger}, in particular (10.10), but
also \cite{AchKre, JonNjaThr, Sim1, Sze2}): $\varphi _1, \ldots, \varphi _n$,
$0\leq \varphi _1< \ldots < \varphi _n < 2\pi$, generate a $(n-1,n, d\sigma )$
Laurent qf
if and only if there is a $\eta \in \mathbb{T}:= \{z\in\mathbb{C}:
\abs{z}=1\}$ such that
\begin{equation}\label{x1.11}	
z\Phi_{n-1}(z) + \eta \Phi_{n-1}^*(z) = \prod\limits_{s=1}^{n}(z-e^{i\varphi_s})
\end{equation}
Here and in what follows
$\Phi_n(z) = z^n + \ldots$ always denotes the monic polynomial of degree $n$ orthogonal with
respect to the normalized measure $d \sigma$, i.e.,
\begin{equation}
\int\limits_{0}^{2\pi}e^{-ik\varphi}\Phi_n(e^{i\varphi})d\sigma(\varphi) = 0 \hspace{1cm} k=0,\ldots,n-1
\end{equation}
It is well known that the $\Phi_n$'s satisfy a recurrence relation of the form
\begin{equation}\label{1.2}
\Phi_n(z) = z\Phi_{n-1}(z) - a_{n-1}\Phi_{n-1}^*(z),
\end{equation}
with $\abs{a_j}<1$ for $j=0,1,2,\ldots$ and where $\Phi_{n-1}^*(z) = z^{n-1}\overline \Phi_{n-1}(\frac{1}{z})$; hence
\begin{equation}\label{1.3}
\Phi_n^*(z) = \Phi_{n-1}^*(z) - \overline a_{n-1} z \Phi_{n-1}(z),
\end{equation}
Furthermore let us recall that for any $\kappa \in\mathbb{T} $

\begin{equation}
z\Phi_{n-1}+\kappa \Phi_{n-1}^* = \prod_{s=1}^{n}(z-e^{i\psi_s})
\label{eq-Nja}
\end{equation}
where $0\leq \psi_1<\ldots<\psi_n<2\pi$, , since $\Phi_{n-1}$ has all zeros in
$\mathbb{D}=\{\abs{z}<1\}$.
Following \cite{JonNjaThr} a polynomial of the form \eqref{eq-Nja} is called
para-orthogonal polynomial.

The Szeg\H{o} quadrature weights are given by, $z_s=e^{i\psi_s}$,
\begin{equation}\label{x1.7}
\mu_s= \frac{1}{2z_s}\frac{(z\Psi_{n-1} -\eta\Psi_{n-1}^*)(z_s)}{(z\Phi_{n-1}+\eta\Phi_{n-1}^*)'(z_s)} > 0
\end{equation}
where $\Psi_n$ is the so-called polynomial of second kind, i.e.
\begin{equation}\label{x1.6}
\Psi_n(z) = \frac{1}{2\pi}\int\limits_{0}^{2\pi}\frac{z+e^{i\varphi}}{z-e^{i\varphi}}(\Phi_n(z)-\Phi_n(e^{i\varphi}))d\sigma(\varphi)
\end{equation}
which satisfies the recurrence relation
\begin{equation}\label{1.2c}
\Psi _n(z) = z\Psi _{n-1}(z) + a_{n-1}\Psi _{n-1}^*(z).
\end{equation}
Note that $a_{n-1}$ in \eqref{1.2} is just replaced by $-a_{n-1}.$
Recently Szeg\H{o} quadrature formulas have been studied and applied by several
authors \cite{BraLiSri, DarGonNja, Gra, Gol, KimRei}. For para-orthogonal
polynomials see \cite{CMV, Sim3, Won}.

As mentioned already, here we study positive $(n-1-m, n, d\sigma )$ qf. There
are two possible approaches, one approach comparable to the algebraic polynomial
approach and the other one via orthogonal polynomials on the unit circle (abbreviated
by OPUC). Each of which has its advantages as we shall see. In the next section we
derive some characterizations without using OPUC's. In fact what concerns Theorem \ref{thm},
we even do not know how to obtain it via OPUC's and what concerns the other
properties as the orthogonality condition given in Lemma \ref{lem}, they follow
in a more natural way than with the help of OPUC's, compare Lemma \ref{lemma2}
and Corollary \ref{corf}. In section 3 we give a simple complete description of
positive trigonometric - resp. Laurent qf with the help of OPUC's. In this
problem the powerfullness of the approach via the unit circle turns out nicely.
We note that section 3 and section 4 can be read without having studied section 2.

\section{Basic characterization}
First let us take a look at the degree of exactness. Similarly as in the case of
algebraic polynomials the degree of exactness is connected with the orthogonality
property of the nodes polynomial. For Szeg\H{o} qf see \cite{BC} and for the associated
orthogonal trigonometric polynomials and their five term recurrence relation see \cite{Jac}.

\begin{lemma}\label{lem}
Let $n \in \mathbb N,$ $n = 2(\tilde{n} + \gamma),$ $\gamma \in \{ 0,1/2 \}.$ A qf of
the form \eqref{1} is exact for $t \in {\mathcal T}_{\tilde{n} + 2 \gamma + l},$
$l \in {\mathbb N}_0,$ if and only if the trigonometric nodes polynomial
$T_{\tilde{n},\gamma}(\varphi)$ satisfies the orthogonality conditions
\begin{equation}\label{L1}
    \int_{0}^{2 \pi} t(\varphi) T_{\tilde{n},\gamma}(\varphi) d\sigma(\varphi) =
    0 {\rm \ for \ } t \in {\mathcal T}_{l,\gamma}
\end{equation}
\end{lemma}

\begin{proof}
Case 1) $n = 2 \tilde{n},$ i.e., $\gamma = 0.$

Sufficiency. First let us consider the case $l = 0.$ Let $\xi \in [0,2 \pi)$ be no zero
of $T_{\tilde{n}}$ and arbitrary otherwise. Since the Lagrange trigonometric polynomials
with respect to the nodes $\varphi_1,...,\varphi_{2 \tilde{n}}, \xi$ are polynomials from
${\mathcal T}_{\tilde{n}}$ the interpolation quadrature formula
\begin{equation}\label{xla1}
   \int_{0}^{2 \pi} t(\varphi) d\sigma = \sum\limits_{s=1}^{2 \tilde{n}} \mu_s t(\varphi_s) + \mu_{2 \tilde{n}+1} t(\xi)
\end{equation}
holds for every $t \in {\mathcal T}_{\tilde{n}}$ and is unique. Now any $q \in {\mathcal
T}_{\tilde{n}}$ can be represented in the form
\begin{equation}\label{xla3}
   q(\varphi) = c T_{\tilde{n}}(\varphi) + t(\varphi)
\end{equation}
where $t \in \mathcal{T}_{\tilde{n}}$ is the unique polynomial satisfying the
interpolation conditions $t(\xi)= 0 $ and $t(\varphi_s) = q(\varphi_s)$ for $s =
1,...,2\tilde{n}$ and $c = q(\xi)/ T_{\tilde{n}}(\xi).$ Hence applying
\eqref{xla1} to $q$, taking into consideration the orthogonality property of
$T_{\tilde{n}}$, we obtain
$$ \frac{1}{2 \pi} \int_0^{2\pi} q d\sigma = \frac{1}{2 \pi} \int_0^{2\pi} c T_{\tilde{n}}
d\sigma + \frac{1}{2 \pi} \int_0^{2\pi} t d\sigma = \sum\limits_{s=1}^{2 \tilde{n}}\mu_s
t(\varphi_s) = \sum\limits_{s=1}^{2 \tilde{n}}\mu_s q(\varphi_s)$$

In case of a higher orthogonality property of $T_{\tilde{n}}$we represent
$q \in {\mathcal T}_{\tilde{n} + l}$ in the form $$ q = v T_{\tilde{n}} + t$$
where $t$ is as above and $v \in {\mathcal T}_l$ and the assertion follows as above.

The necessity part follows immediately by applying the quadrature formula to
$t {T}_{\tilde{n}}.$

Case 2) $n = 2 \tilde{n} + 1,$ i.e., $\gamma = 1/2.$

Sufficiency. Any $q \in {\mathcal T}_{\tilde{n}+1+l}$ can be represented in the form
$$ q = v_{l,1/2} T_{\tilde{n},1/2} + t$$ where $t \in {\mathcal T}_{\tilde{n}}$
interpolates $q$ at the $2 \tilde{n} + 1$ zeros of $T_{\tilde{n},1/2}$ and where
$v_{l,1/2} \in {\mathcal T}_{l,1/2}.$ Applying the unique $(\tilde{n}, 2 \tilde{n} + 1,
d\sigma )$ quadrature formula \eqref{xla1} to $q,$ putting $\xi = \varphi_{2 \tilde{n} + 1},$
the assertion follows by the orthogonality property.

Necessity is immediate again.
\end{proof}

We note if $l < \tilde{n} - 1$ then the orthogonality property \eqref{L1} does not imply
that $T_{\tilde{n},\gamma }$ has all zeros simple and in $[0,2\pi)$ there may
appear also double or complex zeros.

To find a criterion on the nodes polynomial such that the quadrature weights are
positive is more involved (in contrast to the algebraic polynomial case, where the
corresponding counterpart can be derived in one line, see e.g. \cite[(1.3)]{PehMC2}
\begin{thm}
	\label{thm} Let $n\in \N \setminus \{1\}$, $n = 2 (\tilde{n}+\gamma),$ $\gamma \in \{ 0, 1/2 \},$
     and let the qf \eqref{1} be a
    $(\tilde{n}+2\gamma, 2 (\tilde{n}+\gamma), d\sigma )$ qf.
	Denote by $T(\varphi)$
	the trigonometric nodes polynomial, which is from ${\mathcal
	T}_{\tilde{n},\gamma},$ and put
\begin{equation}\label{thm1tilde2}
     S(\psi) = \frac{1}{2\pi} \int_{0}^{2 \pi} \cot
     \left( \frac{\varphi - \psi}{2} \right) (T(\psi) - T(\varphi)) d\sigma(\varphi).
\end{equation}
Then $S(\psi) \in {\mathcal T}_{\tilde{n},\gamma}$ and
\begin{equation}\label{thm1tilde0}
     \mu_s = - \frac{S(\varphi_s)}{2 T^{'}(\varphi_s)} \;\; \text{for}\;\; s = 1,...,n.
\end{equation}
In particular, all quadrature weights $\mu_s$ are positive if and only if $S(\varphi)$
has $n$ simple zeros $\theta_j,$ $0 \leq \theta_1 < \theta_2 < ... < \theta_n < 2 \pi,$
which separate the zeros $\varphi_j$ of $T(\varphi)$ such that $\varphi_1 < \theta_1
< \varphi_2 < ...$ if $sgn(ST)(0^+)> 0$, respectively, $\theta_1 < \varphi_1 < \theta_2 < ...$
if $sgn(ST)(0^+)<0.$

Furthermore, setting
$$\tau_n(e^{i \varphi}) = e^{i \frac{n}{2} \varphi} T(\varphi) \ {\rm and \ }
\omega_n(e^{i \varphi}) = i e^{i \frac{n}{2} \varphi} S(\varphi),$$
the ratio of the polynomials $\omega_n(z)$ and $\tau_n(z)$ has an expansion of the form
\begin{equation}
	-\frac{\omega_n(z)}{\tau_n(z)} = \sum\limits_{s=1}^{n} \mu_s
    \frac{ e^{i \varphi_s} + z}{e^{i \varphi_s} - z} =
    1 + 2\sum\limits_{k=1}^{\tilde{n} + 2 \gamma + l} c_k z^k + ...,
\label{eq-p1}
\end{equation}
where $c_k = \int_{0}^{2 \pi} e^{-i k \varphi} d\sigma(\varphi)$
and $\tilde{n} + 2 \gamma + l$ is the precise degree of exactness of the
qf \eqref{1}.
\end{thm}

\begin{proof}
Let $R(\varphi) \in {\mathcal T}_{\tilde{n},\gamma}$ be the polynomial such that, $z = e^{i \varphi},$
\begin{equation}\label{proofthm1tilde1}
     \frac{i e^{i \frac{n}{2} \varphi} R(\varphi)}{e^{i \frac{n}{2} \varphi} T(\varphi)} =
     \sum\limits_{s=1}^{n} \mu_s \frac{z + e^{i \varphi_s}}{z - e^{i \varphi_s}}
\end{equation}

Putting $\omega_n(e^{i \varphi}) = i e^{i \frac{n}{2} \varphi} R(\varphi)$ and
$\tau_n(e^{i \varphi}) = e^{i \frac{n}{2} \varphi} T(\varphi)$ we may rewrite equation
\eqref{proofthm1tilde1} in the form
\begin{equation}\label{proofthm1tilde-2}
     \frac{\omega_n(z)}{\tau_n(z)} - 1 =
     - 2\sum\limits_{s=1}^{n} \frac{\mu_s}{1 - e^{-i \varphi_s}z}
\end{equation}
where we used the fact that
$$ \sum\limits_{s = 1}^{n} \mu_s = \frac{1}{2 \pi} \int d \sigma = 1$$
since the qf is exact for constants and $\sigma$ is normalized.
Note that $\omega_n(0) = -\tau_n(0)$ or in other words that

\begin{equation}
     T(\varphi) =  a_n \cos \frac{n}{2} \varphi + b_n \sin
     \frac{n}{2} \varphi + ... \;\; \text{and}\;\; R(\varphi) = -b_n\cos \frac{n}{2} \varphi +
     a_n \sin \frac{n}{2} \varphi + ...
\label{eq-p2}
\end{equation}
Furthermore it follows from \eqref{proofthm1tilde-2} that
$$ \frac{\omega_n(e^{i \varphi_s})}{\frac{d \tau_n}{dz}
(e^{i \varphi_s})} = 2 \mu_s e^{i \varphi_s} $$
which yields that
\begin{equation}
	\frac{i R(\varphi_s)}{\frac{dT}{d\varphi}(\varphi_s)} = \frac{e^{-i \frac{n}{2}
	\varphi_s} \omega_n(e^{i \varphi_s})} { \frac{d}{d \varphi} ( e^{-i \frac{n}{2}
	\varphi_s} \tau_n(e^{i \varphi}))_{\varphi = \varphi_s} } = \frac{ \omega_n(e^{i
	\varphi_s}) } { \frac{d}{d \varphi} ( \tau_n(e^{i \varphi}))_{\varphi =
	\varphi_s} } = \frac{2 \mu_s}{i}.
\label{eq-p3}
\end{equation}

The goal is to show that
\begin{equation}\label{proofthm1tilde-1}
      S \equiv  R.
\end{equation}

Let us prove \eqref{proofthm1tilde-1} when $n$ is even.

Using the facts that
\begin{equation*}
    \cot \left( \frac{\varphi - \psi}{2} \right) (\cos \psi - \cos \varphi) =
    \sin \psi + \sin \varphi
\end{equation*}
and that
\begin{equation}
	\frac{\cos k \psi - \cos k \varphi }{\cos \psi - \cos \varphi} = 2 \cos (k-1)
	\psi + \ldots
	\label{eq-p4}
\end{equation}
is a cosine polynomial of degree $k - 1$ in $\psi$ with leading coefficient $2,$
it follows by straightforward calculation that for $k \in \mathbb N$

\begin{equation}
\begin{split}
	&  \int_{0}^{2\pi} \cot \left( \frac{\varphi - \psi}{2} \right) \left(\cos k
	   \psi - \cos k \varphi \right) d \sigma(\varphi) \\
	& = \int_{0}^{2\pi} \cot \left( \frac{\varphi - \psi}{2} \right) (\cos \psi - \cos \varphi)
	  \left(\frac{\cos k \psi - \cos k \varphi }{\cos \psi - \cos \varphi}\right)d
	  \sigma(\varphi)  \\
	&  = (\int_{0}^{2\pi}d\sigma ) \sin k \psi + t(\psi),
\end{split}
\label{eq-p5}
\end{equation}
where $t \in {\mathcal T}_{k-1}$. Furthermore,

\begin{equation}
\begin{split}
	& \int_{0}^{2\pi} \cot \left( \frac{\varphi - \psi}{2} \right) (\sin k \psi -
		\sin k \varphi )d \sigma(\varphi) \\
	& = \sin \psi  \int_{0}^{2\pi} \cot \left( \frac{\varphi - \psi}{2} \right)
	    \left(\frac{\sin k \psi }{\sin \psi } - \frac{\sin k \varphi}{\sin \varphi}\right)
        d \sigma(\varphi) \\
	& + \int_{0}^{2\pi} \cot \left( \frac{\varphi - \psi}{2} \right)(\sin \psi -
	    \sin \varphi )\frac{\sin k \varphi}{\sin \varphi} d\sigma(\varphi) \\
	& = -(\int_{0}^{2\pi}d\sigma ) \cos k \psi + t(\psi)
\end{split}
\label{eq-p6}
\end{equation}
where $t \in {\mathcal T}_{k-1}$, and where in the last equality
we used the facts that $\sin k x / \sin x = 2 \cos (k-1) x +...$ is a cosine polynomial
of degree $k-1$ in $x$ and thus \eqref{eq-p5} can be applied to the first integral and
that the second integral is from ${\mathcal T}_1$. Thus $S(\psi) \in
{\mathcal T}_{\tilde{n}}$ and by \eqref{eq-p2} it follows that $S$ and $R$ have the same
leading coefficients, \ie $R - S \in {\mathcal T}_{\tilde{n} -1}$.
Writing $T(\varphi)$ in terms of $\sin \left(\frac{\varphi - \varphi_s}{2} \right)$
it follows that the integrand from \eqref{thm1tilde2} at $\psi=\varphi _s$ is from
${\mathcal T}_{\tilde{n}}$ and thus we may apply
the quadrature formula to the integrand and obtain
\begin{equation}\label{alt}
    S(\varphi_s) = - \mu_s 2 T^{'}(\varphi_s) = R(\varphi_s), \ s = 1, ..., 2\tilde{n},
\end{equation}
where the last equality follows by \eqref{eq-p3}. Hence $S\equiv R$ and the
assertion is proved for $n$ even.


For odd $n$ we consider
$$ \tilde{S}(\psi) = \frac{1}{2\pi} \int_{0}^{2 \pi} \cot
\left( \frac{\varphi - \psi}{2} \right) ( (q T)(\psi) - (qT)(\varphi) ) d\sigma(\varphi)  $$
where $q$ is any polynomial from ${\mathcal T}_{0,1/2}$ with $q(\varphi) > 0$ on $[0,2\pi).$
Since $qT \in {\mathcal T}_{\tilde{n}+1}$ it follows by \eqref{eq-p5} and \eqref{eq-p6}
that $\tilde{S}\in {\mathcal T}_{\tilde{n}+1}$. Furthermore, we may apply the quadrature
formula and obtain that
\begin{equation}\label{x0}
   \tilde{S}(\varphi_s) = - \mu_s 2 q(\varphi_s) T^{'}(\varphi_s) = q(\varphi_s)
   R(\varphi_s) \;\;\text{for}\;\; s=1,\ldots, 2\tilde{n}+1.
\end{equation}
Next we claim that
\begin{equation}\label{x1}
     \tilde{S}(\psi) = q(\psi) S(\psi)
\end{equation}
To prove \eqref{x1} it suffices to show that
$$ \int_{0}^{2\pi} \cot{\left( \frac{\varphi-\psi}{2} \right)}
\left( q(\psi) - q(\varphi) \right) T(\varphi) d\sigma(\varphi) = 0$$
which follows immediately by $q(\psi) - q(\varphi) = const \sin \left(
\frac{\varphi-\psi}{2} \right)$ and the orthogonality property of $T$ with
respect to ${\mathcal T}_{0,1/2}.$

Therefore, by \eqref{x1} and \eqref{x0}, if
we are able to show that the leading coefficients of $\tilde{S}(\varphi)$ and
$(qR)(\varphi)$ coincide, \ie $\tilde{S} - qR \in \mathcal{T}_{\tilde{n}}$ representation
\eqref{thm1tilde2} follows for $n$ odd also.

Because of  \eqref{eq-p2} we know the leading coefficients of $(qT)(\varphi)$ and
$(qR)(\varphi )$. By applying \eqref{eq-p5} and \eqref{eq-p6} to $(qT)(\psi)-(qT)(\varphi)$ it follows
that the leading coefficients of $\tilde{S}$ are equal to that ones
of $qR$. Hence $\tilde{S} = qR$ and the theorem is proved.


The statement on the positivity of the $\mu_s$ follows immediately by \eqref{alt}.

The last equality in \eqref{eq-p1} follows by partial fraction expansion and \eqref{1.5}.

\end{proof}

Theorem \ref{thm} (which turns out to be crucial in Szeg\H{o}-Kronrod quadrature
\cite{PehSK} for instance) does not show how to obtain a constructive
description of positive trigonometric qf. Such a description will be obtained by
OPUC's in the next section.

\section{The Characterization Theorem}


%

First let us recall that the monic polynomials $\Phi_n(z)$ of degree $n,$
$n \in \mathbb N _0,$ are the polynomials orthogonal with respect to $d \sigma,$
where $\frac{1}{2 \pi} \int_{0}^{2 \pi} d \sigma = 1,$ and that the $a_j$'s are their
recurrence coefficients, see \eqref{1.2}.

\begin{thm}\label{thm1}
The following four statements are equivalent:
\begin{itemize}
\item[a)]
     \begin{equation}\label{1.8}
         \frac{1}{2\pi}\int\limits_{0}^{2\pi}e^{-ik\varphi}d\sigma(\varphi) =
         \sum_{s=1}^{n}\mu_se^{-ik\varphi_s} \hspace{1cm} {\rm for \ } k=0,\ldots,\pm(n-m-1)
     \end{equation}
     where $0\leq \varphi_1<\varphi_2<\ldots<\varphi_n<2\pi$ and $\mu_1,\ldots,\mu_n\in\R^+$.

\item[b)]
     There exists a polynomial $\tilde\Phi_{n-1}(z) = z^{n-1} + \ldots $
     generated by the recurrence relation
     \begin{equation}\label{1.9}
         \tilde\Phi_{j}(z) = z\tilde\Phi_{j-1}(z) - \tilde a_{j-1} \tilde\Phi_{j-1}^*(z),
     \end{equation}
     where
     \begin{equation}\label{1.10}	
         \abs{\tilde a_j}<1 \;\; {\rm for \ } j=0,\ldots,n-2 \;\;{\rm and \ }\;
         \tilde a_j=a_j \;\;	{\rm for \ } j=0,\ldots,n-m-2,
     \end{equation}
     and a $\eta\in\mathbb{T}$ such that
     \begin{equation}\label{1.11}	
         z\tilde\Phi_{n-1} + \eta \tilde\Phi_{n-1}^* = \prod\limits_{s=1}^{n}(z-e^{i\varphi_s})
     \end{equation}

\item[c)]
     There exists a polynomial $q_m(z)=z^m+\ldots$ which has all its zeros in $\mathbb{D}$
     and a $\eta \in \mathbb{T}$ such that
     \begin{equation}\label{1.12}
         z q_m \Phi_{n-m-1} + \eta q_m^*\Phi_{n-m-1}^* = \prod_{s=1}^{n}(z-e^{i\varphi_s})
     \end{equation}

\item[d)]
     There exists a polynomial $p_{n-1}(z)=z^{n-1}+\ldots$ which has all
     its zeros in $\mathbb{D}$ and a $\eta \in \mathbb{T}$ such that
     \begin{equation}\label{1.13}
         z p_{n-1} + \eta p_{n-1}^* = \prod_{s=1}^{n}(z-e^{i\varphi_s})
     \end{equation}
     and
     \begin{equation}\label{1.14}
         \frac{\eta p_{n-1}^* - z p_{n-1}}{\eta p_{n-1}^* + z p_{n-1}} = 1 + 2\sum_{k=1}^{n-m-1}c_k z^k + O(z^{n-m})
     \end{equation}
\end{itemize}
where $c_k = \frac{1}{2\pi}\int\limits_{0}^{2\pi}e^{-ik\varphi}d\sigma(\varphi)$
\end{thm}
\begin{proof}
a)$\Rightarrow$ b) Since $\mu_1,\ldots,\mu_n > 0$ it follows that
the sequence
\begin{equation}
     \tilde c_k =\sum_{s=1}^{n}\mu_s e^{-ik\varphi_s} \hspace{1cm} k=0,\ldots,n-1
\end{equation}
is positive definite, see \cite{AchKre}. Thus the polynomials $\tilde\Phi_j$,
$j=0,\ldots,n-1$, orthogonal with respect to the positive definite sequence $\tilde
c_0,\ldots,\tilde c_{n-1}$ satisfy a recurrence relation of the form (\ref{1.9})
with $\abs{\tilde a_j}<1$ for $j=0,\ldots,n-2$. Since by assumption
\begin{equation}
     \tilde c_k = \frac{1}{2\pi}\int e^{-ik\varphi}d\sigma(\varphi) \hspace{1cm} k=0,\ldots,n-m-1
\end{equation}
it follows that, see e.g. \cite[(3.2)]{Ger}
\begin{equation}
     \tilde a_k=a_k \hspace{1cm} {\rm for\ } k=0,\ldots,n-m-2
\end{equation}
Finally we know \cite[Thm.6.1]{Ger} that for any $a_{n-1}$ with $\abs{a_{n-1}} = 1$
\begin{eqnarray}\label{1.18}
     -\frac{z\tilde\Psi_{n-1}-a_{n-1}\tilde\Psi_{n-1}^*}{z\tilde\Phi_{n-1}+a_{n-1}\tilde\Phi_{n-1}^*} =
     1 + 2\sum_{k=1}^{n-1}\tilde c_k z^k + O(z^n)\nonumber\\[-1.5ex]\\[-1.5ex]
     = - \sum_{s=1}^{n}\mu_s\frac{z+e^{i\varphi_s}}{z-e^{i\varphi_s}} + O(z^n)\nonumber
\end{eqnarray}
Putting $a_{n-1} = -e^{i\varphi_1}\tilde\Phi_{n-1}(e^{i\varphi_1})/\tilde\Phi_{n-1}^*(e^{i\varphi_1})$, that is,
the denominator in (\ref{1.18}) vanishes at $e^{i\varphi_1}$, it follows that,
see \cite[(9.5)]{Ger}
\begin{equation}
     z\tilde\Phi_{n-1}+a_{n-1}\tilde\Phi_{n-1}^* = \prod_{s=1}^n(z-e^{i\varphi_s}).
\end{equation}

b) $\Rightarrow$ a) By (\ref{1.9}) and the first condition of (\ref{1.10})
$\tilde\Phi_{n-1}(z)$ is orthogonal with respect to some positive measure
$\tilde \sigma$. Thus (\ref{1.11}) holds for any $\eta$ with $\abs{\eta}=1$ and
\begin{equation}\label{1.20}
   \frac{1}{2\pi}\int\limits_{-\pi}^{+\pi}e^{-ik\varphi} d\tilde\sigma(\varphi) =
   \sum_{s=1}^n\mu_s e^{-ik\varphi_s} \hspace{1cm}{\rm for\ } k=0,\ldots,n-1
\end{equation}
where $\mu_1,\ldots,\mu_n\in\R^+$. Since $\tilde a_j=a_j$ for $j=0,\ldots,n-m-2$
it follows that $\tilde \Phi_j=\Phi_j$ and $\tilde \Psi_j = \Psi_j$ for $j=0,\ldots,n-m-1$
and thus the first $n-m-1$ moments of $\sigma$ and $\tilde \sigma$ coincide, that is,
\begin{equation}\label{1.21}
   \int e^{-ik\varphi} d \tilde\sigma(\varphi) = \int e^{-ik\varphi}
   d\sigma(\varphi) \hspace{1cm} k=0,\ldots, n-m-1
\end{equation}
which gives by (\ref{1.20}) the assertion.

b) $\Rightarrow$ c) Let
\begin{equation}\label{1.22}	
    q_j(z) = zq_{j-1}(z) - \eta \overline{\tilde a}_{n-1-j}q_{j-1}^*(z)
    \hspace{0.5cm} j=1,2,\ldots,m
\end{equation}
with $q_0(z)=1$ then it follows by induction arguments using (\ref{1.2}) and (\ref{1.3}) that
\begin{equation}\label{1.23}
    z \tilde\Phi_{n-1} + \eta \tilde\Phi_{n-1}^* = z q_j \tilde \Phi_{n-j} + \eta q_j^*\tilde\Phi_{n-j}^*
\end{equation}
In view of \eqref{1.22} $q_j$
has all zeros in $\mathbb{D}$, see e.g. \cite{Sim1}.

We note that it follows analogously by induction arguments that
\begin{equation}\label{1.24}	
    z \tilde\Psi_{n-1} - \eta \tilde\Psi_{n-1}^* = z q_j\Psi_{n-j} - \eta q_j^*\Psi_{n-j}^*
\end{equation}

c) $\Rightarrow$ b) Since $q_m(z)$ has all zeros in $\mathbb{D}$ $q_m$ can be generated
by a recurrence relation of the form
\begin{equation}\label{1.25}	
    q_j(z) = z q_{j-1}(z) - \eta \overline{b}_{j-1}q_{j-1}(z) \hspace{0.5cm} j=1,\ldots,m
\end{equation}
where $\abs{b_j}<1$ for $j=0,\ldots,m-1$ and $q_0(z) = 1$. Hence defining the polynomial
$\tilde \Phi_{n-1}$ by the recurrence relation (\ref{1.2}) with recurrence coefficients
for OPUC's $a_0,\ldots,a_{n-m-1},b_{m-1},\ldots,b_0$ it follows as above that (\ref{1.11}) holds which
implies b).	

c) $\Rightarrow$ d) We put
\begin{eqnarray}\label{1.27}	
    z q_m \Phi_{n-m-1} + \eta q_m^*\Phi_{n-m-1}^* =
    z p_{n-1} + \eta p_{n-1}^* \nonumber \\[-1.5ex]	\\[-1.5ex]
    z q_m \Psi_{n-m-1} - \eta q_m^*\Psi_{n-m-1}^* = z p_{n-1} - \eta p_{n-1}^*	\nonumber
\end{eqnarray}
that is,
\begin{eqnarray}\label{1.28}	
    2 z p_{n-1} = z q_m (\Phi_{n-m-1} + \Psi_{n-m-1}) + \eta q_m^*(\Phi_{n-m-1}^* -
    \Psi_{n-m-1}^*)\nonumber \\[-1.5ex]\\[-1.5ex]
    2 \eta p_{n-1}^* = z q_m (\Phi_{n-m-1} - \Psi_{n-m-1}) + \eta
    q_m^*(\Phi_{n-m-1}^* + \Psi_{n-m-1}^*) \nonumber
\end{eqnarray}
Let us demonstrate that $p_{n-1}^*$ has no zero in $\overline{\mathbb{D}}$. It
is known, see (1.3.82), (3.2.62) and (3.2.63) from \cite{Sim1} that $\Phi_{n-m-1}^* +
\Psi_{n-m-1}^*$ has no zeros on $\overline{\mathbb{D}}$. We claim that
\begin{equation}\label{1.29}	
   \abs{\frac{z(\Phi_{n-m-1}-\Psi_{n-m-1})}{\Phi_{n-m-1}^* + \Psi_{n-m-1}^*}}<1 	
   \hspace{0.5cm} {\rm for \ } \abs{z}\leq 1.	
\end{equation}
By the maximum principle it suffices to show (\ref{1.29}) for $\abs{z}=1$. Since on
$\abs{z}=1$
\begin{equation}\label{1.30}	
    \abs{\frac{z(\Phi_{n-m-1}-\Psi_{n-m-1})}{\Phi_{n-m-1}^* + \Psi_{n-m-1}^*}} = 	
    \abs{\frac{\Phi_{n-m-1}^* - \Psi_{n-m-1}^*}{\Phi_{n-m-1}^* + \Psi_{n-m-1}^*}}	
\end{equation}
and since by \cite[(10.3)]{Ger}
\begin{equation}	
{\rm Re} \frac{\Psi_{n-m-1}^*(e^{i\varphi})}{\Phi_{n-m-1}^*(e^{i\varphi})}>0
\end{equation}
it follows that the last expression in (\ref{1.30}) is less than $1$ and thus
\eqref{1.29} is proved. Since $q_m$ has all zeros in $\mathbb{D}$ there holds
$|q_m(z)/q_m^*(z)|<1$ on $\overline{\mathbb{D}}$, which implies by \eqref{1.29}
and \eqref{1.28} that $p_{n-1}^*$ has no zero in $\overline{\mathbb{D}}$.
Finally, relation (\ref{1.14}) follows with the help of the relation, see
\cite[Prop. 3.2.9]{Sim1}
\begin{equation}	
    \Phi_{n-m-1}(z)F(z) + \Psi_{n-m-1}(z) = O(z^{n-m-1})
\end{equation}
and
\begin{equation}	
    \Phi_{n-m-1}^*(z)F(z) - \Psi_{n-m-1}^*(z) = O(z^{n-m})
\end{equation}
and (\ref{1.27}).

d) $\Rightarrow$ c) The proof runs as in \cite[p. 939]{Peh3}. To make the paper
selfcontained we sketch it. According to \cite[Theorem 18.2]{Ger}) there exists
a function $h:\C\rightarrow\C$ which is analytic in $\mathbb{D}$ and satisfies the
inequality $\abs{h(z)}<1$ for $z\in \mathbb{D}$, such that
\begin{equation}\label{0.34}
    \frac{\eta p_{n-1}^*(z)-z p_{n-1}(z)}{\eta p_{n-1}^*(z)+z p_{n-1}(z)} = -\frac{z
    \Psi_{n-1-m}(z)h(z)-\Psi_{n-m-1}^*(z)}{z \Phi_{n-1-m}(z)h(z)+\Phi_{n-m-1}^*(z)}.
\end{equation}
Isolating $h$ from this equality we find that $h$ can be represented by
\begin{equation}
   h(z) = \frac{q_l(z)}{q_l^*(z)}, \hspace{1cm} {\rm where \ } q_l(z) =
   \prod_{i=1}^{l}(z-z_i), \ \ z_i\in \mathbb{D},
\end{equation}
$q_l$ has real coefficients, and $l\geq m$.

Let us assume that $l>m$. Then it follows from (\ref{0.34}) that
$-z \Psi_{n-1-m}h + \Psi_{n-m-1}^*$ and $z \Phi_{n-1-m}h + \Phi_{n-m-1}^*$ have
$(l-m)$ common zeros on $\abs{z}=1$, which implies that
\begin{equation}
    \frac{\Phi_{n-m-1}}{\Phi_{n-m-1}^*}+\frac{\Psi_{n-m-1}}{\Psi_{n-m-1}^*}
\end{equation}
has $(l-m)$ zeros on $\abs{z}=1$. But this is impossible, since (see \cite[p. 4]{Ger})
\begin{equation}
   \Omega_{n-m-1}\Psi_{n-m-1}^* + \Psi_{n-m-1}\Phi_{n-m-1}^* = K z^{n-1-m},
\end{equation}
where $K\in\R^+$.

a) $\Rightarrow$ d) Let
\begin{equation}\label{0.24}
    c_k= \frac{1}{2\pi}\int\limits_{0}^{2\pi} e^{-ik\varphi}d\sigma(\varphi) =
    \sum_{s=1}^{n}\mu_s e^{-ik\varphi_s} \hspace{1cm} {\rm for \ } k=0,\ldots,n-m-1.
\end{equation}
Hence, by expansion at $z=0$,
\begin{equation}\label{0.25}
    - \sum_{s=1}^{n}\mu_s \frac{z+e^{i\varphi_s}}{z-e^{i\varphi_s}} =
    1 + 2 \sum_{k = 1}^{n-1-m} c_k z^k + O(z^{n-m}).
\end{equation}
Now the left hand side is a rational function with real part zero on $\abs{z}=1$ whose numerator and denominator is a
polynomial of degree $n$. Thus it may be written as
\begin{equation}\label{0.26}
    \frac{\eta p_n^*(z)-p_n(z)}{\eta p_n^*(z)+p_n(z)} = - \sum_{s=1}^{n}\mu_s
    \frac{z+e^{i\varphi_s}}{z-e^{i\varphi_s}} = 1 + \ldots,
\end{equation}
hence
\begin{equation}\label{0.27}
    p_n(z) = zp_{n-1}(z)
\end{equation}
Moreover by (\ref{0.26}), partial fraction expansion and the suppositions on $\mu_s$, using also
$z\frac{d}{dz}f(z) = -i\frac{d}{d\varphi}f(e^{i\varphi})$, we obtain that
\begin{equation}\label{0.28}
    \mu_s=\frac{(\eta p_n^*-p_n)(e^{i\varphi_s})}{2i \frac{d}{d\varphi}(p_n +
    \eta p_n^*)(e^{i\varphi})_{\varphi=\varphi_s}}>0.
\end{equation}
Thus the trigonometric polynomials ${\rm Im}
\{\eta ^{-1/2}e^{-i\frac{n}{2}\varphi}p_n(e^{i\varphi})\}$ and \\
${\rm
Re}\{\eta ^{-1/2}e^{-i\frac{n}{2}\varphi}p_n(e^{i\varphi})\}$ have strictly interlacing
zeros, which implies by the argument principle that $p_n$ has all zeros in
$\mathbb{D}$. By (\ref{0.26}) and (\ref{0.25}) the implication is proved.

d) $\Rightarrow$ a) Since $p_n:=zp_{n-1}$ has all zeros in $\mathbb{D}$ it follows that
$\eta p_n^*\pm p_n$ has all zeros on $\partial \mathbb{D}$ and that the zeros of
$\eta p_n^* - p_n$ and $\eta p_n^* + p_n$ strictly interlace. Hence \eqref{0.28} and thus
(\ref{0.26}) holds with $\mu_s>0$ for $s=1,\ldots,n$. Expanding the second expression from
(\ref{0.26}) in a series at $z=0$ part a) follows.
\end{proof}

Obviously Theorem \ref{thm1} could be proved in a shorter way, that is, by
showing the implications a)$\Rightarrow$b)$\Rightarrow \ldots
\Rightarrow$d)$\Rightarrow$a) only. But by proving all equivalences separately
we get additional insights, for instance:

\begin{remark}
In the proof of Theorem \ref{thm1}, c)$\Rightarrow$b) we have shown that
$\tilde\Phi_{n-1}$ is generated by the recurrence coefficients
$a_0,\ldots,a_{n-m-1}, b_{m-1}, \ldots, b_0$ when $q_m$ is generated by
\begin{equation}
    q_j(z) = z q_{j-1}(z) - \eta \bar b_{j-1}q_{j-1}^*(z) \hspace{1cm} j=1,\ldots,m
\end{equation}
$q_0(z)=1$, with $|b_j|<1$, $j=0,\ldots,m-1$. Note that every polynomial $q_m$ which has all
zeros in $\mathbb{D}$ can be generated by a recurrence relation of the form \eqref{1.2}, see
\cite[Thm. 9.2]{Ger}.
\end{remark}
The ideas of proof for Theorem \ref{thm1}c) and d) go back essentially to the
author's papers \cite{Peh2,Peh3} where positive qf on $[-1,1]$ have been characterized by
transforming the problem to the unit circle. For the connection with qf on $[-1,1]$
see Section 5 below.

Finally we point out that recurrence relation \eqref{1.9}, under condition \eqref{1.10},
offers a simple, constructive way to generate positive trigonometric $(n - 1 - m, n, d \sigma)$ qf.

With the help of Theorem \ref{thm1} we get easily some information on the
location of the nodes of positive qf on the unit circle.

\begin{cor} \label{cor} Let $l, m\in \mathbb{N}_0, m\leq l$, and let
$z\Phi_{n-1-l}+\kappa \Phi _{n-1-l}^* = \prod _{\nu =1}^{n-1-l} (z-e^{i\psi _\nu})$,
$\kappa \in \mathbb{T}$, be a Szeg\H{o} polynomial. Suppose that
$T_n(z) = \prod_{s=1}^n (z-e^{i\varphi _s})$ generates a positive $(n-1-m, n, d\sigma )$
qf. Then in every interval $(\psi _\nu, \psi _{\nu +1})$, $\nu = 0,\ldots,l$, there
is at least one $\varphi _s$.

Furthermore, every common zero of $T_n$ and $z\Phi _{n-1-l}+\kappa \Phi_{n-1-l}^* $
is a zero of $\kappa q_l(z)-\eta q_l^*(z)$, where $q_l$ and $\eta$ is given by the
representation $T_n = zq_l\Phi_{n-1-l} + \eta q_l^*\Phi_{n-1-l}^*$.
\end{cor}

\begin{proof} By Theorem \ref{thm1} and the obvious fact that every $(n-1-m,
n, d \sigma )$ qf is a $(n-1-l, n, d\sigma )$ qf for $l\geq m$, $T_n$ has a
representation of the form
\begin{equation}
	T_n = zq_l\Phi_{n-1-l} + \eta q_l^*\Phi_{n-1-l}^*
	\label{eq-g1}
\end{equation}
where $q_l$ has all zeros in $\mathbb{D}$.
Hence $T_n(e^{i\varphi _s}) = 0$ if and only if for some $k\in \mathbb{Z}$

\begin{equation}
	 \arg \eta + (2k+1)\pi = \arg \frac{zq_l\Phi_{n-1-l}}{q_l^*\Phi_{n-1-l}^*}(e^{i\varphi _s}) =
	 \arg \frac{q_l}{q_l^*}(e^{i\varphi _s}) + \arg \frac{z\Phi_{n-1-l}}{\Phi_{n-1-l}^*}(e^{i\varphi _s})
	\label{eq-g0}
\end{equation}
and the condition for a zero of the Szeg\H{o}-polynomial looks similarily.
Now recall that the argument of a Blaschke product is strictly increasing on $|z|=1$
with respect to $\varphi $. Thus on $[\psi _\nu, \psi _{\nu +1}]$
$\arg \frac{z\Phi_{n-1-l}}{\Phi_{n-1-l}^*}$ increases by
$2\pi$ and, taking a look at the RHS of \eqref{eq-g0}, $\arg \frac{zq_l\Phi_{n-1-l}}
{q_l^*\Phi_{n-1-l}^*}$ increases by more than $2\pi$ on $[\psi _\nu, \psi _{\nu +1}]$,
that is, $T_n$ has at least one zero in $(\psi _\nu, \psi _{\nu +1})$ and the assertion is proved.
\end{proof}

For the special case that $T_n=z\Phi _{n-1} + \eta \Phi _{n-1}^*$ the zero
property given in the above Corollary was proved by Simon \cite[Theorem
2.3]{Sim3} by a completely different approach. For related statements in the real
case see \cite{BD, PehMC, Peh4, Sze1}.

%

For the convenience of the reader let us reformulate parts of Theorem \ref{thm1} in
terms of trigonometric polynomials.
\begin{cor}
$\varphi_1,...,\varphi_n,$ $0 \leq \varphi_1 < \varphi_2, ... < \varphi_n < 2 \pi,$
are the nodes of a positive trigonometric $(n - 1 - m, n, d\sigma)$ quadrature formula
if and only if there exists a polynomial $\tilde{\Phi}_{n-1}(z)$ satisfying
\eqref{1.9} and \eqref{1.10} and a $\eta \in \mathbb T$ such that, $z = e^{i \varphi},$
\begin{equation}\label{polx}
   \eta^{-1/2} z^{-n/2} \prod\limits_{s=1}^{n} (z -  e^{i \varphi_s}) =
   2 Re \{ z^{-\frac{n}{2}} \eta^{1/2} \tilde{\Phi}^{*}_{n-1}(z) \}
\end{equation}
\end{cor}

\begin{proof}
The statement follows immediately by Theorem \ref{thm1} a) and b) and the obvious
relations, $z = e^{i \varphi},$
\begin{equation}
\begin{split}
  & \eta^{-1/2} z^{- \frac{n}{2}} \left[ z \tilde{\Phi}_{n-1}(z) + \eta \tilde{\Phi}^{*}_{n-1}(z) \right] =
    2 Re \{ z^{-\frac{n}{2} + 1} \eta^{-1/2} \tilde{\Phi}_{n-1}(z) \} \\
  & = 2 Re \{ z^{-\frac{n}{2} } \eta^{1/2} \tilde{\Phi}^{*}_{n-1}(z) \} \} .\\
\end{split}
\end{equation}
\end{proof}

We need the following lemma which can be extracted easily from
\cite{Sze2} or \cite[pp. 144-145]{Ger2}.

\begin{lemma}\label{lemma2}
Let $\mu$ be a positive measure on $[0, 2 \pi],$ let $\tilde{n} \in \mathbb N,$
$\gamma \in \{0,1/2 \}$ and $n = 2(\tilde{n} + \gamma).$ Let
$P_{n-1}(z) = z^{n-1} + ...$ be such that
\begin{equation}\label{tilde1}
   \int_0^{2 \pi} e^{- i k \varphi} P_{n-1}(e^{i \varphi}) d\mu(\varphi) = 0 {\rm \ for \ } k = 0,...,n-2
\end{equation}
Then for any $\eta \in {\mathbb T}$
\begin{equation}\label{tilde2}
   \int_0^{2 \pi} t(\varphi) Re \{ \eta^{1/2} e^{- i \frac{n}{2}\varphi} P^*_{n-1}(e^{i
   \varphi}) \} d\mu(\varphi) = 0 {\rm \ for \ } t \in \mathcal
   T_{\tilde{n}-1,\gamma}
\end{equation}
\end{lemma}

\begin{proof}
Relation \eqref{tilde2} is equivalent to
\begin{equation}
\begin{split}
    & \int_0^{2 \pi} e^{- i \left( \frac{n}{2} -1 \mp (k + \gamma) \right)  \varphi} P_{n-1}(e^{i \varphi}) d\mu(\varphi) \\
  + & \eta \int_0^{2 \pi} e^{- i \left( \frac{n}{2} \mp (k + \gamma) \right)  \varphi} P_{n-1}^{*}(e^{i \varphi}) d\mu(\varphi)
      = 0 {\rm \ for \ } k = 0,..., \tilde{n} - 1.\\
\end{split}
\end{equation}
Now the first integral is zero by \eqref{tilde1} and thus the second
integral also by taking the complex conjugate of it.
\end{proof}


With the help of the above Lemma we get an explicit weight function depending on $n$
with respect to which the trigonometric nodes polynomial has maximal orthogonality.
Furthermore we obtain the orthogonality property \eqref{L1} of
the nodes polynomial \eqref{polx} by an approach via OPUC's.

\begin{cor} \label{corf} Let $n = 2(\tilde{n} + \gamma),$ where $\tilde{n} \in {\mathbb N}$ and
$\gamma \in \{ 0, 1/2 \}$. If $\varphi_1,...,\varphi_n,$ $0 \leq \varphi_1 <
\varphi_2, ... < \varphi_n < 2 \pi,$ are the nodes of a positive trigonometric
$(n-1-m,n,d\sigma)$ qf then there exists a polynomial $q_m(z)$ which has all
zeros in $|z| < 1$ and a $\eta \in \mathbb T$ such that $Re \{ \eta^{1/2} z^{-
\frac{n}{2}} ( q_m^{*}\Phi^{*}_{n-1-m} )(z) \},$ $z=e^{i \varphi},$ is the nodes
polynomial and
\begin{equation}\label{xto1}
\begin{split}
     & \int_0^{2 \pi} t(\varphi) Re \{ \eta^{1/2} e^{- i \frac{n}{2} \varphi}
       (q_m^{*} \Phi^{*}_{n-1-m})(e^{i \varphi})  \}
       \frac{d\varphi}{|(q_m^{*} \Phi^{*}_{n-1-m})(e^{i \varphi})|^2} \\
     & = 0 \ {\rm \ for \ } t \in \mathcal T_{\tilde{n}-1, \gamma}.
\end{split}
\end{equation}
Moreover,
\begin{equation}\label{xto2}
     \int_0^{2 \pi} t(\varphi) Re \{ \eta^{1/2} e^{- i \frac{n}{2} \varphi }
     (q_m^{*} \Phi^{*}_{n-1-m})(e^{i \varphi})  \} d\sigma(\varphi)
     = 0 \ {\rm \ for \ } t \in \mathcal T_{\tilde{n}-1-m,\gamma}
\end{equation}
\end{cor}

\begin{proof}
By Theorem \ref{thm1} $q_m \Phi_{n-1-m}$ has all zeros in $\mathbb{D}$, hence, see e.g. \cite{Ger, Sze1},
$$
     \int_0^{2 \pi} e^{- i k \varphi} (q_m \Phi_{n-1-m})(e^{i \varphi})
     \frac{d\varphi} {|(q^*_m \Phi^*_{n-1-m})(e^{i \varphi})|^2} = 0 \ {\rm \ for \ } k =
     0,...,n-2.
$$
Applying Lemma \ref{lemma2} the orthogonality property \eqref{xto1} follows.

Since $|q^*_m(e^{i \varphi})|^2$ is from ${\mathcal T}_m$ we obtain from \eqref{xto1},
by putting $t(\varphi) = \tilde{t}(\varphi) |q^*_m(e^{i \varphi})|^2,$
$\tilde{t} \in {\mathcal T}_{\tilde{n} - 1 - m, \gamma},$ that \eqref{xto2} holds
with respect to the measure $d \varphi / |\Phi^*_{n - 1 - m}(e^{i \varphi})|^2$
instead of $d \sigma.$ But, see \cite[Thm 2.2, p. 198]{Fre},
\begin{equation*}
  \int_{0}^{2 \pi} e^{\pm i k \varphi} d \sigma (\varphi) =
  \int_{0}^{2 \pi} e^{\pm i k \varphi} \frac{d \varphi}{|\Phi^*_{n - 1 - m}(e^{i \varphi})|^2}
  {\rm \ for \ } k = 0, ..., n - 1 - m
\end{equation*}
and thus relation \eqref{xto2} follows, taking into consideration the fact that the orthogonality
property \eqref{xto2} depends on the first $ n - m - 1 = 2\tilde{n} + 2\gamma - m - 1$ trigonometric moments of $d \sigma$ only.

\end{proof}

We note that in general ${\rm Im} \{ z^{-\frac{n}{2}} \eta^{1/2}
\tilde{\Phi}^*_{n-1}(z) \},$ $\tilde{\Phi}_{n-1}$ given by \eqref{1.9} and \eqref{1.10},
does not coincide with ${\rm Im} \{ z^{-\frac{n}{2}} \eta^{1/2} (q_m^{*} \Phi^{*}_{n-1-m})(z) \} $,
since
$$z \tilde{\Phi}_{n-1} - \eta \tilde{\Phi}^{*}_{n-1} = z \tilde{q}_m \Phi_{n-1-m}
- \eta \tilde{q}^{*}_m \Phi^{*}_{n-1-m},$$ where $\tilde{q}_m$ is given by
$\tilde{q}_j(z) = z \tilde{q}_{j-1}(z) +  \eta \bar{\tilde{a}} _{n-1-j} \tilde{q}^{*}_{j-1}(z),$
$j = 1,2,...,m$ with $\tilde{q}_0(z) = 1,$ which should be compared with \eqref{1.22}.


\section{Asymptotics of weights}

As a consequence of Theorem \ref{thm1} we obtain the following representation
of the quadrature weights needed in what follows. Note that only the $\Phi_n$'s and $q_m$'s and no
$\Psi_n$'s appear in the representation.

\begin{cor}
The quadrature weights $\mu_s$ from (\ref{1.8}) can be represented in the form
\begin{equation}
	\mu _s = \frac{ -\eta K_{n-m-1} z_s^{n-1}\abs{q_m(z_s)}^2 }{(z\Phi_{n-m-1}q_m -
	\eta \Phi_{n-m-1}^* q_m^*)(z_s)(z\Phi_{n-m-1}q_m + \eta \Phi_{n-m-1}^* q_m^*)'(z_s)}
	\label{eq-ro30}
\end{equation}
where $z_s = e^{i \varphi_s}$ and $K_{n-m-1} = 2\prod _{j=0}^{n-m-2}(1-|a_j|^2)$
\end{cor}

\begin{proof}
By \eqref{x1.7}, \eqref{1.23} and \eqref{1.24} it follows that
\begin{equation}
	\mu _s =  \frac{1}{2z_s}\frac{ (z \Psi_{n-m-1}
	q_m  - \eta \Psi_{n-m-1}^* q_m^* )(z_s)}{(z\Phi_{n-m-1}q_m + \eta \Phi_{n-m-1}^*
	q_m^*)'(z_s)}
	\label{eq-ro31}
\end{equation}
Now, at the zeros of $z\Phi_{n-m-1}q_m + \eta  \Phi_{n-m-1}^* q_m^*$ we have
\begin{eqnarray*}
& (z\Phi_{n-m-1}q_m - \eta  \Phi_{n-m-1}^* q_m^*)( z \Psi_{n-m-1}
	q_m -\eta  \Psi_{n-m-1}^* q_m^* ) = \\
	& -2\eta zq_mq_m^*(\Phi_{n-m-1}\Psi_{n-m-1}^* + \Psi_{n-m-1}\Phi_{n-m-1}^*)  \\
	& = -2\eta z q_m q_m^* K_{n-m-1}z^{n-m-1}
 \end{eqnarray*}
where in the last equality we used the known fact that \cite[(5.6)]{Ger}
\begin{equation}
	\Phi _{n-m-1}\Psi _{n-m-1}^* + \Psi _{n-m-1}\Phi _{n-m-1}^* = K_{n-m-1}z^{n-m-1}
	\label{eq-nr00}
\end{equation}
which gives by \eqref{eq-ro31} the assertion.
\end{proof}

\begin{thm}
Let $d \sigma(\varphi) = f(\varphi) d \varphi$ be positive and from $Lip \gamma,$
$0 < \gamma \leq 1$ on $[\alpha, \beta] \subseteq [0, 2 \pi].$ Suppose that
$\varphi_{1,n}, ... , \varphi_{n,n},$ $0 \leq \varphi_{1,n} < \varphi_{2,n} <
... < \varphi_{n,n} < 2 \pi$ generates a sequence of positive $(n- m(n)-1,n,w)$
qf with quadrature weights $\mu_{s,n}.$ Furthermore let us assume that the
associated $q_{m(n)}$ satisfy uniformly on $[\alpha, \beta]$ $$\lim \limits_n
\left( 1 - \frac{2}{n} Re \{ e^{i \varphi} \frac{ q^{{*}^{'}}_{m(n)}(e^{i
\varphi}) } { q^{*}_{m(n)}(e^{i \varphi}) } \} \right) = g(\varphi)$$
Then uniformly for $\varphi _{s,n}\in [\alpha + \varepsilon, \beta -
\varepsilon], \varepsilon >0,$
\begin{equation}
    \frac{1}{n \mu_{s,n}} = \frac{g(\varphi_{s,n})}{f(\varphi_{s,n})} + o(1)
\end{equation}
\end{thm}

\begin{proof}
First let us note that \eqref{eq-ro30} can be written in the form
\begin{equation}\label{0.42}
\begin{split}
    & - 4 \mu_s = \\
    & \frac{K_{n-m-1}\abs{q_m(e^{i\varphi_s})}^2}
         {{\rm Im}  \{ \overline {e^{-i\frac{n}{2}\varphi_s}(\eta^{1/2}q_m^*\Phi_{n-m-1}^*)(e^{i\varphi_s})}\}
    \frac{d}{d\varphi}
    {\rm Re}\{e^{-i\frac{n}{2}\varphi}(\eta^{1/2}q_m^*\Phi_{n-m-1}^*)(e^{i\varphi})\}_{\varphi=\varphi _s}}\\
\end{split}
\end{equation}
Now
\begin{equation}\label{0.46}
\begin{split}
    & \ \ \ \ {\rm Im} \{\overline{e^{-i\frac{n}{2}\varphi}(\eta^{1/2}q_m^{*}\Phi_{n-m-1}^*)(e^{i\varphi})}\}
      \frac{d}{d\varphi}{\rm Re} \{e^{-i\frac{n}{2}\varphi}(\eta^{1/2}q_m^{*}\Phi_{n-m-1}^*)
      (e^{i\varphi})\} \\
    & - {\rm Re}\{e^{-i\frac{n}{2}\varphi}(\eta^{1/2}q_m^{*}\Phi_{n-m-1}^*)(e^{i\varphi})\} \frac{d}{d\varphi}
      {\rm Im} \{e^{-i\frac{n}{2}\varphi}(\eta^{1/2}q_m^{*}\Phi_{n-m-1}^*)(e^{i\varphi})\} \\
    & = {\rm Im} \{(\overline{e^{-i\frac{n}{2}\varphi}(\eta^{1/2}q_m^{*}\Phi_{n-m-1}^*)(e^{i\varphi})})
      \frac{d}{d\varphi}(e^{-i\frac{n}{2}\varphi}(\eta^{1/2}q_m^{*}\Phi_{n-m-1}^*)(e^{i\varphi})) \} \\
    & = {\rm Im} \{ - i \frac{n}{2} |(q^{*}_{m} \Phi^{*}_{n - m - 1}  )(e^{i \varphi})|^2  +
      \overline{(q^{*}_m \Phi^{*}_{n - m - 1})(e^{i \varphi})}\frac{d}{d \varphi}
      \left( (q^{*}_m \Phi^{*}_{n - m - 1})(e^{i \varphi}) \right) \}
\end{split}
\end{equation}
and thus, since the $e^{i \varphi_s}$ are the zeros of ${\rm Re}
\{e^{-i\frac{n}{2}\varphi}(\eta^{1/2}q_m^{*}\Phi_{n-m-1}^*)(e^{i\varphi})\}$,
we may replace the denominator in \eqref{0.42} by the last expression from \eqref{0.46}
which yields
\begin{equation}
    \mu_s\abs{(q_m^{*}\Phi_{n-m-1}^*)(e^{i\varphi_s})}^2 =
    \frac{K_{n-m-1}\abs{q_m^*(e^{i\varphi_s})}^2}{4{\rm Im}
    \{i\frac{n}{2}-ie^{i\varphi_s}\frac{(q_m^{*}\Phi_{n-m-1}^*)'
    (e^{i\varphi_s})}{(q_m^{*}\Phi_{n-m-1}^*)(e^{i\varphi_s})}\}}
\end{equation}

Now by the assumptions on $f$ it is known, see e.g. \cite{Sze1, Sim1} that uniformly on
$[\alpha + \varepsilon, \beta - \varepsilon], \ \varepsilon > 0,$
\begin{equation}
	\lim\limits_{n \to \infty} \frac{K_{n-m(n)-1}}{2|\Phi^*_{n-m(n)-1}(e^{i
\varphi})|^2} = \frac{1}{|D(e^{i \varphi},\sigma)|^2} = f(\varphi)
\label{p-1}
\end{equation}
taking into consideration that
\begin{equation}
	\lim _n K_{n-m(n)-1}/2= 1/D(0,\sigma),
	\label{eq-p8}
\end{equation}
where $D(z,\sigma)$ is the so-called Szeg\H{o} function, that is,
\begin{equation}
	D(z) =
	\exp\{\frac{1}{2\pi}\int\limits_{0}^{2\pi}\frac{e^{i\varphi}+z}{e^{i\varphi}-z} \log
	f(\varphi )d \varphi \}
	\label{eq-f2}
\end{equation}

Next let us prove that \eqref{p-1} and \eqref{eq-p8} imply that uniformly on $[\alpha +
\varepsilon,\beta - \varepsilon], \varepsilon>0$,
\begin{equation}\label{eq5.24}
\frac{1}{n}
\abs{\frac{\Phi_{2n-1-m}^{*'}(e^{i\varphi})}{\Phi_{2n-1-m}^{*}(e^{i\varphi})}}
\underset{n \to \infty}{\longrightarrow } 0
\end{equation}
Indeed
by the local version of Bernstein's inequality

\begin{equation}
	\begin{split}
		& \max_{\varphi\in[\alpha + \varepsilon ,\beta - \varepsilon]}|\frac{d}{d\varphi}
          \Phi_{2n-m-1}^{*}(e^{i\varphi})| \leq const \times  \\
        & \qquad \qquad ((2n-m-1)\max_{\varphi\in[\alpha,\beta]}
	      \abs{\Phi_{2n-m-1}^{*}(e^{i\varphi})-\Phi_{[\sqrt{2n-m-1 }]}^{*}(e^{i\varphi})} \\
        & \qquad \qquad + \sqrt{2n-m-1}\max_{\varphi\in[\alpha,\beta]}\abs{\Phi_{[\sqrt{2n-m-1}]}^{*}(e^{i\varphi})})
\end{split}
\label{eq-new1}
\end{equation}
which gives by
\eqref{p-1} and \eqref{eq-p8}, in conjunction with the facts that
$D(e^{i\varphi})\not=0$ on $[\alpha,\beta]$ and $m(n)\leq n$, relation
\eqref{eq5.24}.
\end{proof}

\section{Connection to positive qf on $[-1,1]$}

\begin{remark} Let $m,n\in\mathbb{N}_0, 0\leq m\leq 2n-1$. Obviously, using the
symmetry with respect to $\pi$ we have:
\begin{equation}
    \int\limits_{-1}^{+1}p(x)d\psi(x) = \sum_{s=1}^n \lambda_s p(x_s)
    \hspace{1cm} {\rm for \ }p \in \mathbb{P}_{2n-1-m},
\end{equation}
where $-1<x_n<\ldots<x_1<1$ and $\lambda_1,\ldots,\lambda_n \in \mathbb{R}^+$,
if and only if
\begin{equation}
   \int\limits_{-\pi}^{+\pi}t(\varphi)d\sigma(\varphi) = \sum_{s=1}^n \lambda_s
   (t(\varphi_s)+t(-\varphi_s)) \hspace{1cm} {\rm for \ }t \in \mathcal{T}_{2n-1-m},
\end{equation}
where the $\varphi_s=\arccos{x_s}$ satisfy $0<\varphi_1<\varphi_2<\ldots<\varphi_n<\pi,$
and $\sigma$ is given by
$$
\sigma(\varphi) =
\left\{
\begin{aligned}
   \psi(1) - \psi(\cos \varphi) & {\rm \ \ for \ } & \ 0 \leq \varphi \leq \pi,\\
   \psi(\cos \varphi) - \psi(1) & {\rm \ \ for \ } & - \pi \leq \varphi \leq 0.\\
\end{aligned}
\right.
$$
\end{remark}

Note that, $x=\frac{1}{2}(z+\frac{1}{z})$, $$2^n z^n \prod\limits_{s=1}^n(x-x_s) =
\prod\limits_{s=1}^{n}(z-e^{i\varphi_s})(z-e^{-i\varphi_s})$$.

By this equivalence and Theorem \ref{thm1} we obtain immediately a complete
description of positive qf of degree of exactness $2n-1-m$ on $[-1,1]$. Indeed,
Theorem \ref{thm1}c) and d) will lead to characterizations given by the author
in \cite{Peh2}. In fact in \cite{Peh2} we transformed the problem in the
proof to the unit circle solved it there, that is, proved Theorem \ref{thm1}c)
and d) in the real case and transformed it back. It can be shown that the
characterization by the recurrence relation, that is, by Theorem \ref{thm1}b) yields,
using the connection between the recurrence coefficients of OPUC's and the recurrence
coefficients of polynomials orthogonal on $[-1,1]$ (see \cite[(13.1.7)]{Sim2} or
\cite[Thm. 31.1]{Ger}), the characterization of positive qf on $[-1,1]$ by the
three term recurrence relation given in \cite{Peh3,Peh4} by the author, see also \cite{PehMC2}.

When we consider Radau or Lobatto qf on $[-1,1]$, see \cite[Section 4]{PehMC2} that is,
that one or two nodes
are the boundary points $\pm 1$ the corresponding equivalence is as follows: let
us consider the Radau case, that is, the boundary point $1$ is a node - analogously
the other cases are obtained - then
\begin{equation}
    \int p(x) d\psi(x) = \sum_{s=1}^n \lambda_s p(x_s) +
    \lambda_0 p(+1) \hspace{1cm} {\rm for \ }p\in \mathbb{P}_{2n-m},
\end{equation}
where $-1<x_n<\ldots<x_1<1$ and $\lambda_0,\lambda_1,\ldots,\lambda_n \in \mathbb{R}^+,$ if and only if
\begin{equation}
    \int\limits_{-\pi}^{+\pi}t(\varphi)d\sigma(\varphi) = \frac{\lambda_0}{2}t(\varphi_0) +
    \sum_{s=1}^n \lambda_s (t(\varphi_s)+t(-\varphi_s)) \hspace{1cm} {\rm for \ }t \in \mathcal{T}_{2n-m},
\end{equation}
where the $\varphi_s=\arccos{x_s}$ satisfy $\varphi_0=0<\varphi_1<\varphi_2<\ldots<\varphi_n<\pi$
and where $\lambda_0,\lambda_1,\ldots,\lambda_n \in \mathbb{R}^+$.



\begin{thebibliography}{99}
\bibitem{AchKre} N.I. Achieser and M. Krein, {\em The L-problem of moments},
In: Achieser, N.I., Krein, M. (eds.) Some Questions in the Theory of Moments.
Translation of Mathematical Monographs, Vol. 2, Amer. Math. Soc., Providence, RI, 1962


\bibitem{BD} A.F. Beardon and K.A. Driver, {\em The zeros of linear combinations
of orthogonal polynomials}, J. Approx. Theory 137 (2005), no.2, 179 - 186

\bibitem{BC} R. Cruz-Barroso, P. Gonz\'alez Vera and O. Nj\r astad, {\em On
bi-orthogonal systems of trigonometric functions and quadrature formulas for
periodic integrals}, Numer. Alg. 44 (2007), 309-333.

\bibitem{BraLiSri} C.F. Bracciali, Xin Li and A. Sri Ranga, {\em Real orthogonal
polynomials in frequency analysis}, Math. Comp. 74(249), 341-362 (2004)


\bibitem{CMV} M.J. Cantero, L. Moral and L. Vel\'azquez, {\em Measures and
para-orthogonal polynomials on the unit circle}, East J. Approx. 8 (2002), 447-464



\bibitem{DarGonNja} L. Daruis, P. Gonz\'alez-Vera and O. Nj\r astad, {\em Szeg\H{o}
quadrature formulas for certain Jacobi type weight functions}, Math. Comp. 71, 683-701 (2002)

\bibitem{Fre} G. Freud, {\em Orthogonal Polynomials}, Pergamon Press, Oxford (1971)

\bibitem{Ger} Ya.L. Geronimus, {\em Polynomials orthogonal on a circle and their applications},
Amer. Math. Soc. Transl. 3, 1-78 (1962)

\bibitem{Ger2} Ya.L. Geronimus, {\em Orthogonal Polynomials}, Consultants
Bureau, New York, 1961

\bibitem{Gol} L. Golinskii, {\em Quadrature Formula and zeros of para-orthogonal
polynomials on the unit circle}, Acta Math. Hungar. 96(3), 169-186 (2002)

\bibitem{Gra} W.B. Gragg, {\em Positive definite Toeplitz matrices, the Arnoldi
process for isometric operators, and Gaussian quadrature on the unit circle}(in Russian),
In: Nicholaev, E.S. (ed.) Numerical Methods in Linear Algebra, pp. 16-32. Moscow
University Press, Moscow (1982). Published in English in slightly modified form
in J. Comput. Appl. Math. 46, 183-198 (1993)





\bibitem{Jac} D. Jackson, {\em Series of orthogonal polynomials and orthogonal
trigonometric sums}, Annals of Mathematics (2) 34 (1933), 527 - 545, 799 - 814


\bibitem{JonNjaThr} W.B. Jones, O. Nj\r astad and W.J. Thron, {\em Moment theory,
orthogonal polynomials, quadrature and continued fractions associated with the unit
circle}, Bull. London Math. Soc. 21, 113-152 (1989)

\bibitem{KimRei} Sun-Mi Kim and L. Reichel, {\em Anti-Szeg\H{o} quadrature rules},
Math. Comp. 76(258), 795-810 (2007)


\bibitem{MilCveSta} G.V. Milovanovi\'c, A.S. Cvetkovi\'c and M.P. Stani\'c, {\em
Trigonometric orthogonal systems and quadrature formulae},
Comput. Math. Appl. 56 (2008), no. 11, 2915-2931

\bibitem{Peh2} F. Peherstorfer, {\em Characterization of positive quadrature
formulas}, SIAM J. Math. Anal. 12, 935-942 (1981)

\bibitem{Peh3} F. Peherstorfer, {\em Characterization of quadrature formulas II},
SIAM J. Math. Anal. 15, 1021-1030 (1984)

\bibitem{PehMC} F. Peherstorfer, {\em Linear combinations of orthogonal
polynomials generating positive quadrature formulas}, Math. Comp. 55(1990), 231-241


%
%
%

\bibitem{Peh4} F. Peherstorfer,  {\em On positive quadrature
formulas}, ISNM, 112, Birk\"{a}user, Basel, 297-313 (1993)



\bibitem{PehMC2} F. Peherstorfer, {\em Positive quadrature formulas III:
Asymptotics of weights}, Math. Comp. 77 (2008), 2241-2259


\bibitem{PehSK} F. Peherstorfer, {\em Szeg\H{o}-Kronrod quadrature on the unit
circle}, manuscript


\bibitem{Sim1} B. Simon, {\em Orthogonal polynomials on the unit
circle, Part 1: Classical theory}, AMS Colloquium Series, American mathematical
Society, Providence, RI, 2005

\bibitem{Sim2} B. Simon, {\em Orthogonal polynomials on the unit circle, Part 2:
Spectral theory}, AMS Colloquium Series, American mathematical Society, Providence, RI, 2005

\bibitem{Sim3} B. Simon, {\em Rank one perturbations and the zeros of
paraorthogonal polynomials on the unit circle}, J. Math. Anal. Appl. 329 (2007), 376-382


\bibitem{Sze1} G. Szeg\H{o}, {\em Orthogonal Polynomials}, 4th ed., Amer. Math. Soc.,
Providence, 1975

\bibitem{Sze2} G. Szeg\H{o}, {\em On bi-orthogonal systems of trigonometric polynomials},
Magyar Tud. Akad. Mat. Kutat 8(1964), 255-273


\bibitem{Won} M.-W.L. Wong, {\em First and second kind paraorthogonal polynomials and their zeros},
J. Approx. Theory 146 (2007), 282 - 293
\end{thebibliography}
\end{document}